\documentclass[11pt]{article}
\usepackage[cp1251]{inputenc}
\usepackage[T2A]{fontenc}
\usepackage[english]{babel}
\usepackage{blindtext}
\usepackage[margin=1cm,left=1.5cm,includefoot]{geometry}
\usepackage[mathcal]{eucal}
\usepackage{amsfonts}
\usepackage{amssymb}
\usepackage{ amssymb }
\usepackage{amsmath}
\usepackage{hyperref}
\usepackage{abstract}
\usepackage[symbol]{footmisc}

\usepackage{amsthm}
\usepackage{amscd}
\usepackage{mathrsfs}

\theoremstyle{plain}

\newtheorem{lem}{Lemma}[section]
\newtheorem{foll}{Corollary}[section]

\newtheorem{theor}{Theorem}
\newtheorem{conj}{Conjecture}

\newenvironment{solve}{\begin{proof}[Proof]}{\end{proof}}

\begin{document}

\renewcommand{\thefootnote}{\fnsymbol{footnote}}

\begin{center}
\textsc{\Large Diophantine properties of fixed points of Minkowski question mark function.}\\ 
\,
\,
\large Dmitry Gayfulin\footnote[1]{Research is supported by RNF grant No. 19-11-00001.}\footnote[2]{The author is a Young Russian Mathematics award winner and would like to thank its sponsors and jury.}, Nikita A. Shulga\footnote[3]{This research was supported by RFBR grant 18-01-00886.}
\end{center}

\begin{abstract}
\noindent
We consider irrational fixed points of the Minkowski question mark function $?(x)$, that is
irrational solutions of the equation $?(x)=x$. It is easy to see that there exist at least two such points. Although 
it is not known if there are other fixed points, we prove that the smallest and the greatest fixed points have irrationality measure exponent equal to 2. We give more precise results about the approximation properties of these fixed points. Moreover, in Appendix we introduce a condition from which it follows that there are only two irrational fixed points.
\end{abstract}

\section{Introduction}
For $x\in [0,1]$ we consider its continued fraction expansion
$$x=[a_1,a_2,\ldots,a_n,\ldots]=  \cfrac{1}{a_1+\cfrac{1}{a_2+\cdots}}, \,\,\, a_j \in \mathbb{Z}_+$$
which is unique and infinite when $x\not\in\mathbb{Q}$  and finite for rational $x$. 
Each rational $x$ has just two representations
$$
x=[a_1,a_2,\ldots,a_{n-1}, a_n],  
\,\,\,\,\,\text{and}\,\,\,\,\,
x=[a_1,a_2,\ldots,a_{n-1}, a_n-1,1]
,
\,\,\,\, \text{where}\,\,\,\, a_n \ge 2
.
$$
By 
$$\frac{p_k}{q_k}:=[a_1,\ldots,a_k]$$ 
we denote the $k$th convergent fraction to $x$.
By $B_n$ we denote the $n$th level of the Stern-Brocot tree, that is $$B_n:=\{ x=[a_1,\ldots,a_k]:a_1+\ldots+a_k=n+1\}.$$
In \cite{Mink} Minkowski introduced a special  function $?(x)$ which may be defined as the limit  distribution function
of sets $B_n$. This function was rediscovered and studied by many authors (see \cite{kinney},\cite{kesse},\cite{alkauskas},\cite{dushistova},\cite{paradis}).
For irrational $x=[a_1,a_2,\ldots,a_n,\ldots]$   the formula
\begin{equation} \label{eq:1}
?(x)=\sum\limits_{k=1}^{\infty} \frac{(-1)^{k+1}}{2^{a_1+\ldots+a_k-1}}
\end{equation}
introduced by Denjoy \cite{Den,Den1} and Salem \cite{Salem}
may be considered as one of the equivalent definitions of the function $?(x)$. If $x$ is rational, then the infinite series in (\ref{eq:1}) is replaced by a finite sum. Note that $?([0;a_1,\ldots,a_t+1])=?([0;a_1,\ldots,a_t,1])$ and hence $?(x)$ is well-defined for rational numbers too. It is known that $?(x)$ is a continuous 
strictly increasing function, its derivative $?'(x)$  exists almost everywhere in $[0,1]$ in the sense of Lebesgue measure, and $?'(x)=0 $
for $ x\in \mathbb{Q}$.

As $?(x)$ is continuous, $?'(0)=?'(\frac{1}{2})=0$ and
$$?(0) = 0, \,\,\,\,?\left(\frac{1}{2}\right)  =\frac{1}{2},\,\,\,\, ?(1) = 1,$$
 we see that there exist at least two points $ x_1 \in  \left(0,\frac{1}{2}\right)$
 and 
 $ x_2 \in  \left(\frac{1}{2},1\right)$ with
 $$
 ?(x_j) = x_j,\,\,\,\, j = 1,2.
 $$
 A folklore conjecture states that
\begin{conj}
\label{mainconj}
The  Minkowski question mark function $?(x)$ has exactly five fixed points. There is only one irrational fixed point of $?(x)$ in the interval $(0, \frac{1}{2})$.
\end{conj}
This conjecture has not yet been proved (for certain discussion  see survey  preprint by Moshchevitin \cite{Moshchevitin}). Our computations show that if there is more then one fixed point in the interval $(0; \frac{1}{2})$, then the first 5400 partial quotients in the continued fraction expansion of these numbers coincide. 
Although we do not know if there are exactly two irrational fixed points of $?(x)$, we are able to say something about Diophantine properties of some of them. In the present paper we give explicit lower bounds for the irrationality measure of the smallest and the greatest fixed points from $(0,\frac{1}{2})$ interval, that is lower bounds of the form
$$
\Bigl|x - \frac{p}{q}\Bigl|>\frac{1}{q^2\cdot I(q)}=\frac{1}{q^{2+\delta(q)}},   \,\,\,\,\, \delta(q)\ge
0$$
satisfied by all $p,q\in\mathbb{Z},q\ge q_0$, where the dependence $I(q)$ on $q$ is explicit and $q_0$ is given. Usually the infimum  $\inf_{q\in \mathbb{Z}_+} (2+\delta(q))$ is called irrationality measure (or exponent) of $x$. 

As $?(1-x)=1-?(x)\ \forall x\in[0,1]$, the set of fixed point of the Minkowski question mark function is symmetric with respect to the point $\frac{1}{2}$. Therefore, one can study the fixed points on $[0,\frac{1}{2}]$ interval only.
\section{Main results}
Our first result establishes some properties of the continued fraction expansion of certain fixed points of $?(x)$. \\
\begin{theor}
\label{th1}
Let $x=[a_1,\ldots,a_n,\ldots]$ be the smallest or the greatest fixed point of Minkowski question mark function on the interval $(0,\frac{1}{2})$. Then $a_1=2$ and
\begin{equation}
\label{eq:2}
\,\,\, a_{n+1}\le\sum\limits_{i=1}^{n}a_i
\end{equation}
for all $ n\in\mathbb{N}$.
\end{theor}
We give a proof of Theorem \ref{th1} in Section \ref{th1proof}. The following theorem is a stronger version of Theorem \ref{th1}. It uses some new geometrical considerations.
\begin{theor}
\label{th1improved}
Denote $\kappa_1=2\log_2(\frac{\sqrt{5}+1}{2})-1\approx 0.38848383\ldots$.  Let $x$ be fixed point from Theorem \ref{th1}, then 
\begin{equation}
\label{eq:2impr}  
\,\,\, a_{n+1}<\kappa_1\sum\limits_{i=1}^{n}a_i+2\log_2\biggl(\sum\limits_{i=1}^{n}a_i\biggr)
\end{equation}
for all $ n\ge 1$.
\end{theor}
Formula \eqref{eq:2impr} gives an explicit irrationality measure lower estimate for the fixed points under considerations.
\begin{theor}
\label{th2}
Let $x$ be fixed point from Theorem \ref{th1}, then $\exists q_0\in\mathbb{Z}_{+}$ such that
\begin{equation*}
\biggl|x-\frac{p}{q}\biggr|>\frac{1}{\biggl(\frac{2\kappa_1}{\log2}\log{q}+O(\log\log{q})\biggr)q^2}
\end{equation*}
for all  $q>q_0\in\mathbb{N}, p\in\mathbb{N}$. 
\end{theor}
We give a proof of Theorem \ref{th2} in Section \ref{th2proof}. The $O(\log\log{q})$ term is calculated explicitly there.

The following statement reduces the problem of fixed points of $?(x)$ to the properties of values of $?(x)$ at rational points only.
\begin{theor}
\label{theorappendix}
Conjecture \ref{mainconj} follows from the inequality
\begin{equation}
\label{ge2q2}
\biggl|?\biggl(\frac{p}{q}\biggr)-\frac{p}{q}\biggr|>\frac{1}{2q^2}
\end{equation}
for all $p,q\in\mathbb{Z}_{+}$ with $q\ge q_0$ for some $q_0\in\mathbb{Z}_{+}$.
\end{theor}
We prove Theorem \ref{theorappendix} in Appendix.

\section{Preliminaries}
\label{prelim}
By $F_k$ we denote $k$th Fibonacci number, that is 
$$F_0=F_1=1 , F_{n+1}=F_{n}+F_{n-1}.$$
\begin{lem}
\label{ratval}
Let $?([a_1,a_2,\ldots,a_{n-1}])=[b_1,b_2,\ldots,b_k]$ and $\sum\limits_{i=1}^{n-1}a_i=s+1>2$, then $\sum\limits_{i=1}^{k}b_i>s+1$.
\end{lem}
\begin{solve}
Reducing \eqref{eq:1} to a common denominator we get
\begin{equation}
\label{eq:3}
?([a_1,a_2,\ldots,a_{n-1}])=\frac{2^{a_2+\ldots+a_{n-1}}-2^{a_3+\ldots+a_{n-1}}+\ldots+(-1)^{n-1}\cdot 2^{a_{n-1}}+(-1)^{n}}{2^{\sum\limits_{i=1}^{n-1}a_i-1}}.
\end{equation}
Here the denominator is equal to $2^s$, and the numerator is an odd number. Let us consider the level $B_s$ of the Stern-Brocot tree, which contains the number $[a_1,a_2,\ldots,a_{n-1}]$. The greatest denominator on this level is equal to $F_{s+2}$. We know that for all $s>2$ one has $F_{s+2}<2^s$. This means that the image of the number
$[a_1,a_2,\ldots,a_{n-1}]$, given by the formula \eqref{eq:3}, belongs to level $B_{s+k}$ for some $k\ge1$, since the denominator of the image is greater than the greatest denominator on the level $B_s$.
\end{solve}
\noindent
Letting $S([a_1,\ldots,a_n])=a_1+\ldots+a_n$, Lemma \ref{ratval} is actually proving the inequality $S(?(x))>S(x)$ for every $x\in\mathbb{Q}\cap(0,1),\ x\neq\frac{1}{2}$.
\begin{foll}
The Minkowski question mark function has exactly 3 rational fixed points: $0,\,\, \frac{1}{2}$ and $1.$
\end{foll}
\begin{solve}
We see that $F_{s+2}=2^s$ only for $s=0,1$, that is for numbers from the $0th$ and the $1st$ levels of the Stern-Brocot tree, which only contain the numbers $0, \frac{1}{2}$ and $1$. For every other rational number, the  sum of its partial quotients increases under the map $?(x)$. So the number is not mapped onto itself.
\end{solve}
The following lemma about the values of Minkowski function at rational points is related to a famous statement known as "Folding lemma" (see \cite{shallit}).
\begin{lem}
\label{fold}
Let $s$ be an arbitrary nonnegative integer and
$$
?([a_1,a_2,\ldots,a_{n-1}])=[b_1,b_2,\ldots,b_k],\,\,\,b_k\neq1.
$$
Consider the number 
$$
\theta =[a_1,a_2,\ldots,a_{n-1}, a_n],\,\,\,\text{where}\,\,\,a_n = \sum\limits_{i=1}^{n-1}a_i+s, \,\,\, s\ge0.$$
Then
\begin{enumerate}
\item{If $n\equiv k(mod$  $2)$, then $?(\theta)=[b_1,b_2,\ldots,b_{k-1},b_k-1,1,2^{s+1}-1,b_k,\ldots,b_1].$}
\item{If $n\equiv k+1(mod$  $2)$ then $?(\theta)=[b_1,b_2,\ldots,b_k,2^{s+1}-1,1,b_k-1,b_{k-1},\ldots,b_1].$}
\end{enumerate}
\end{lem}
\begin{solve}
We know that $b_k\neq1$. Let us choose one of the representations $\frac{p_l}{q_l}=[b_1,b_2,\ldots,b_k]$ or $\frac{p_l}{q_l}=[b_1,b_2,\ldots,b_k-1,1]$ so that the length  $l$ of the continued fraction
expansion
is of the same parity as $n+1$, that is 
$l\equiv n+1\pmod{2}, $ and $ l = k $ or $ l = k+1$. 
From \eqref{eq:3}  we know that $q_l=2^{\sum\limits_{i=1}^{n-1}a_i-1}$.\\
Without loss of generality suppose  that $n\equiv k+1\pmod{2}$, then
$$?(\theta)=?([a_1,a_2,\ldots,a_{n-1}])+\frac{(-1)^{n+1}}{2^{2\cdot \sum\limits_{i=1}^{n-1}a_i-1+s}}=\frac{p_l}{q_l}+\frac{(-1)^{n+1}}{2^{s+1}\cdot q_l^2}=\frac{p_l}{q_l}+\frac{(-1)^{l}}{2^{s+1}\cdot q_l^2}=\frac{p_lq_l2^{s+1}-(p_lq_{l-1}-q_lp_{l-1})}{2^{s+1}q_l^2}=$$
$$=\frac{p_l(2^{s+1}-\frac{q_{l-1}}{q_l})+p_{l-1}}{q_l(2^{s+1}-\frac{q_{l-1}}{q_l})+q_{l-1}}=[b_1,b_2,\ldots,b_k,2^{s+1}-\frac{q_{l-1}}{q_l}]=[b_1,b_2,\ldots,b_k,2^{s+1}-1,1,b_k-1,b_{k-1},\ldots,b_1].$$
In the last equality we use $1-[b_k,\ldots,b_1]=[1,b_k-1,b_{k-1},\ldots,b_1]$.
\end{solve}
\noindent
\textbf{Remark 1.} H.Niederreiter  \cite{Niederreiter} proved that if $m$ is a power of $2$, then there exists an odd integer $a$ with $1\le a\le m$ such that all partial quotients in continued fraction expansion of $\frac{a}{m}$ are bounded by $3$. In fact, he took iterations of the Minkowski question mark function of the continued fractions of a special form, where each partial quotient is equal to the sum of all previous ones or to the sum of all previous ones plus $1$. 

\begin{lem}
\label{nextpar}
Let  $a_1,\ldots,a_{n-1}$ be the partial quotients of a fixed point $x$, then, depending on the parity of $n$, the next partial quotient $a_n$ satisfies one of the following systems
\begin{enumerate}
\item{If $n$ is even, then $a_n$ satisfies $\begin{cases} [a_1,\ldots,a_{n-1},a_n]<?([a_1,\ldots,a_{n-1},a_n+1]), \\ [a_1,\ldots,a_{n-1},a_n+1]>?([a_1,\ldots,a_{n-1},a_n]). \end{cases}$}
\item{If $n$ is odd, then $a_n$ satisfies $\begin{cases} [a_1,\ldots,a_{n-1},a_n]>?([a_1,\ldots,a_{n-1},a_n+1]), \\ [a_1,\ldots,a_{n-1},a_n+1]<?([a_1,\ldots,a_{n-1},a_n]). \end{cases}$\\  }
\end{enumerate}
\end{lem}
\begin{solve}
Let $[a_1,\ldots,a_n,\ldots]$ be the  fixed point. Let us show 1).\\
Since $n$ is even, from the continued fraction theory we know
$$[a_1,\ldots,a_n]<[a_1,\ldots,a_n,\ldots]<[a_1,\ldots,a_n+1].$$
That is $[a_1,\ldots,a_n]$ and $[a_1,\ldots,a_n+1]$ lie on opposite sides with respect to $x$, and hence their images lie on different sides too, and we have
$$\begin{cases} [a_1,\ldots,a_{n-1},a_n]<?([a_1,\ldots,a_{n-1},a_n+1]) ,\\ [a_1,\ldots,a_{n-1},a_n+1]>?([a_1,\ldots,a_{n-1},a_n]). \end{cases}$$
Case 2) can be treated similarly, since for odd $n$ one has 
$$
[a_1,\ldots,a_n+1]<[a_1,\ldots,a_n,\ldots]<[a_1,\ldots,a_n].
$$
\end{solve}
The following lemma localizes fixed points.
\begin{lem}
\label{2537}
All fixed points of ?(x) inside the interval $(0,\frac{1}{2})$ belong to  the interval $(\frac{2}{5},\frac{3}{7})$.
\end{lem}
\begin{solve}
First of all, we will show that there are no fixed points on the interval $(0,\frac{1}{3})$. Decompose the interval $(0,\frac{1}{3})$ into the union of subintervals $(0,\frac{1}{3})\setminus \mathbb{Q}=\bigcup\limits_{n=3}^{\infty}(\frac{1}{n+1},\frac{1}{n}) \setminus \mathbb{Q} $. Assume that for some $n_0$  there exists $ x_0 \in (\frac{1}{n_0+1},\frac{1}{n_0})$ such that $?(x_0)=x_0$. Then
$$ \frac{1}{n_0+1}<x_0=?(x_0)<?\left(\frac{1}{n_0}\right)=\frac{1}{2^{n_0-1}}.$$
 So  $n_0+1>2^{n_0-1}$, 
 and this is not true for
 $\forall n_0 \ge 3$.\\
This means that the first partial quotient is equal to $2$ ($1$ is excluded, since we are on the interval $(0,\frac{1}{2})$).

Now we show that there are no fixed points on the interval $(\frac{4}{9},\frac{1}{2})$. Consider the decomposition $(\frac{4}{9},\frac{1}{2})\setminus \mathbb{Q}=\bigcup\limits_{n=4}^{\infty}(\frac{n}{2n+1},\frac{n+1}{2n+3}) \setminus \mathbb{Q} $. Assume that for some $n_0$  there exists $ x\in (\frac{n_0}{2n_0+1},\frac{n_0+1}{2n_0+3})$ such that $?(x)=x$.
$$
\frac{n_0+1}{2n_0+3}>x_0=?(x_0)>?\left(\frac{n_0}{2n_0+1}\right)=\frac{1}{2}-\frac{1}{2^{1+n_0}}.
$$
We get $2^{n_0}<2n_0+3$, which holds only for $n_0=1,2,3$.

To show that there all fixed points are on the interval $(\frac{2}{5},\frac{3}{7})$, we observe that $?(?(\frac{2}{5}))=?(\frac{3}{8})=\frac{5}{16}<\frac{1}{3}$ and $?(?(\frac{3}{7}))=?(\frac{7}{16})=\frac{29}{64}>\frac{4}{9}$.
\end{solve}
Lemma \ref{2537} means that the continued fraction expansion of every fixed point in $(0,\frac{1}{2})$ is of the form $[2,2,\ldots]$.\\    
The next statement is an obvious property of continuous functions. We formulate it without a proof.
\begin{lem}
\label{contfunc}
Let $f(x)$ be a  continuous function. Consider an interval $[a,b]$ such that there are no fixed points inside $(a,b)$. Then $f(x)-x$ does not change sign on $(a, b)$.
\end{lem}

\section{Proof of Theorem  \ref{th1}}
\label{th1proof}
 Let us prove this theorem for the leftmost fixed point of $?(x)$ on the interval $(0,\frac{1}{2})$, which we denote by $x=[a_1,\ldots,a_n,\ldots]$. Our proof goes by contradiction. Assume that there exists $n\ge3$ such that  
 \begin{equation}\label{main}
 a_n\ge\sum\limits_{i=1}^{n-1}a_i.
 \end{equation}
 Consider $[a_1,\ldots,a_{n-1}]$ and let $?([a_1,\ldots,a_{n-1}])=[b_1,\ldots,b_k],\ b_k\ge 2$. Now we distinguish cases 1) - 4) with several subcases.
 In each of them we will deduce a contradiction.
 We present a very detailed exposition of the case 1) below.
 Cases 2) - 4) are quite similar. We establish them with less details.
 \\
1) $n$ odd, $k$ odd. Then by Lemma \ref{fold}, $?([a_1,a_2,\ldots,a_n])=[b_1,b_2,\ldots,b_{k-1},b_k-1,1,2^{s+1}-1,b_k,\ldots,b_1]$.\\
By Lemma \ref{nextpar} $a_n$ should satisfy \begin{equation} \label{123} \begin{cases} [a_1,\ldots,a_{n-1},a_n]>?([a_1,\ldots,a_{n-1},a_n+1])=[b_1,b_2,\ldots,b_{k-1},b_k-1,1,2^{s+2}-1,b_k,\ldots,b_1], \\ [a_1,\ldots,a_{n-1},a_n+1]<?([a_1,\ldots,a_{n-1},a_n])=[b_1,b_2,\ldots,b_{k-1},b_k-1,1,2^{s+1}-1,b_k,\ldots,b_1]. \end{cases} \end{equation}
Now we consider the possible three subcases 1.1), 1.2) and 1.3).\\
1.1) $k\le n-1$.  Then by Lemma \ref{ratval} we have $\sum\limits_{i=1}^{k}b_i>\sum\limits_{i=1}^{n-1}a_i$, hence there exists $ i\in \{1,\ldots,k\}$ such that $a_i\neq b_i$. By considering a partial quotient with the smallest index $i\le k$ for which $a_i\neq b_i$, we get that the system \eqref{123} is incompatible by the rules of comparison of continued fractions.\\
1.2) $k>n$. If there is $ i\in \{1,\ldots,n-1\}$ such that $a_i\neq b_i$, then similarly to case 1.1) system \eqref{123} is incompatible. Hence we assume that $a_i= b_i$ for all $ i\in \{1,\ldots,n-1\}$. Let us consider the four possible variants for $b_n$.\\
1.2.1) $b_n\le a_n-1$. The first inequality of \eqref{123} can be rewritten as
\begin{equation}\label{1234}[a_1,\ldots,a_{n-1},a_n]>[a_1,a_2,\ldots,a_{n-1},b_n,\ldots,b_{k-1},b_k-1,1,2^{s+2}-1,b_k,\ldots,b_1].\end{equation}
But $n$ is odd and $b_n\le a_n-1$, so \eqref{1234} cannot be true.\\
1.2.2) $b_n\ge a_n+2$. The second inequality of \eqref{123} can be rewritten as
\begin{equation}\label{12345}[a_1,\ldots,a_{n-1},a_n+1]<[a_1,a_2,\ldots,a_{n-1},b_n,\ldots,b_{k-1},b_k-1,1,2^{s+1}-1,b_k,\ldots,b_1].\end{equation}
\\
But $n$ is odd and $b_n\ge a_n+2$, so \eqref{12345} cannot be true.\\
1.2.3) $b_n=a_n+1$. The second inequality in \eqref{123} can be rewritten in this case in the form
$$[a_1,\ldots,a_n+1]<[a_1,\ldots,a_n+1,b_{n+1},\ldots,b_{k-1},b_k-1,1,2^{s+1}-1,b_k,\ldots,b_1]. $$
This inequality fails since the value of continued fraction is always less than  an odd convergent.\\
1.2.4) $b_n=a_n$. Now the equality  $?([a_1,\ldots,a_{n-1}])=[b_1,\ldots,b_k]$   gives
$$?([a_1,\ldots,a_{n-1}])=[a_1,\ldots,a_{n-1},a_n,\ldots,b_k]>[a_1,\ldots,a_{n-1}].$$ 
The last inequality is due to the fact that $n$ is odd.
But $[a_1,\ldots,a_{n-1}]<x $ is also a convergent for our fixed point $x$, so we have
 $$[a_1,\ldots,a_{n-1}]<?([a_1,\ldots,a_{n-1}])<[a_1,\ldots,a_n,\ldots]=x.$$
This contradicts  Lemma \ref{contfunc} because by Lemma \ref{2537} we have $\frac{2}{5}<x<\frac{3}{7}$ and $?(\frac{2}{5})=\frac{3}{8}<\frac{2}{5}$, hence $?(y)<y$ should hold for every $y\in(0,x)$, and in particular for $y=[a_1,\ldots,a_{n-1}]$, which is not possible.\\
1.3) $k=n$. Similarly to the case 1.2) we get $a_i= b_i$ 
for all $ i = 1,\ldots,n-1$.
Let us consider the following four possible subcases.\\
1.3.1) $b_n\le a_n$.  In this case we deduce a contradiction similarly to the case 1.2.1).\\
1.3.2) $b_n\ge a_n+3$. In this case we deduce a contradiction similarly to the case 1.2.2).\\
1.3.3) $b_n=a_n+1$. Equality  $?([a_1,\ldots,a_{n-1}])=[b_1,\ldots,b_k]$ under the assumptions  of this case gives 
 $$?([a_1,\ldots,a_{n-1}])=[a_1,\ldots,a_{n-1},a_n+1]>[a_1,\ldots,a_{n-1}].$$
 But this contradicts Lemma \ref{contfunc} as in the case 1.2.4).\\
1.3.4) $b_n=a_n+2$. The system \eqref{123} gives now in particular
$$[a_1,\ldots,a_n+1]<[a_1,\ldots,a_n+1,1,2^{s+1}-1,a_n,\ldots,a_1].$$
This inequality fails since odd convergents are always greater than the value of the continued fraction.\\\\
2) n even, k even. Then by  Lemma \ref{fold} $?([a_1,a_2,\ldots,a_n])=[b_1,b_2,\ldots,b_{k-1},b_k-1,1,2^{s+1}-1,b_k,\ldots,b_1]$.\\
By  Lemma \ref{nextpar} $a_n$ should satisfy \begin{equation}\label{system10} \begin{cases} [a_1,\ldots,a_{n-1},a_n]<?([a_1,\ldots,a_{n-1},a_n+1])=[b_1,b_2,\ldots,b_{k-1},b_k-1,1,2^{s+2}-1,b_k,\ldots,b_1], \\ [a_1,\ldots,a_{n-1},a_n+1]>?([a_1,\ldots,a_{n-1},a_n])=[b_1,b_2,\ldots,b_{k-1},b_k-1,1,2^{s+1}-1,b_k,\ldots,b_1]. \end{cases} \end{equation}
As before, there are three possible cases.\\
2.1) $k\le n-1$. In this case we deduce a contradiction similarly to the case 1.1).\\
2.2) $k>n$. In the same way as in 1.2), we come to four options:\\
2.2.1) $b_n\le a_n-1$. In this case we deduce a contradiction similarly to the case 1.2.1).\\
2.2.2) $b_n\ge a_n+2$. In this case we deduce a contradiction similarly to the case 1.2.2).\\
2.2.3) $b_n=a_n+1$. The system \eqref{system10} gives now 
$$[a_1,\ldots,a_{n-1},a_n+1]>?([a_1,\ldots,a_{n-1},a_n])=[a_1,\ldots,a_n+1,b_{n+1},\ldots,b_{k-1},b_k-1,1,2^{s+1}-1,b_k,\ldots,b_1].$$
This second inequality fails since even convergents are always smaller than the value of the continued fraction.\\
2.2.4) $b_n=a_n$. Let us consider the even convergent  $[a_1,\ldots,a_n]$ of our fixed point, or rather its image. Taking into account the assumptions and using Lemma \ref{fold}, we obtain\\
$$?([a_1,\ldots,a_n])=[a_1,\ldots,a_n,b_{n+1},\ldots,b_k-1,1,2^{s+1}-1,b_k,\ldots,b_1]>[a_1,\ldots,a_n].$$ We get a contradiction in the same way as in 1.2.4), since we found a number $[a_1,\ldots,a_n]$ less than the {\it the smallest} fixed point, whose image is greater than this number.\\
2.3) $k=n$.  Similarly to the case 1.2) we come to the assumption that $a_i= b_i$ for all $i\in \{1,\ldots,n-1\}$ and so we consider the following four possible subcases:\\
2.3.1) $b_n\le a_n$. In this case we deduce a contradiction similarly to the case 1.2.1).\\
2.3.2) $b_n\ge a_n+3$. In this case we deduce a contradiction similarly to the case 1.2.2).\\
2.3.3) $b_n=a_n+1$. We have 
$$
?([a_1,\ldots,a_n])=[a_1,\ldots,a_n,1,2^{s+1}-1,a_n,\ldots,a_1]>[a_1,\ldots,a_n].
$$
Now we get a contradiction as in 1.2.4).\\
2.3.4) $b_n=a_n+2$. The system \eqref{system10} gives now in particular
$$[a_1,\ldots,a_n+1]>[a_1,\ldots,a_n+1,1,2^{s+1}-1,a_n,\ldots,a_1].$$
This inequality fails since even convergents are always smaller than the value of the continued fraction.\\\\
Cases 
$$
\text{3)}\,\, n \text{ even},\,\, k \text{ odd}
\,\,\,\,\,\,\,\,\,\,\, \text{and}\,\,\,\,\,\,\,\,\,\, 
\text{4)}\,\, n \text{ odd},\,\, k \text{ even}
$$
 follow by a similar argument. For example, in  case 3)
 by Lemma \ref{fold}, $$?([a_1,a_2,\ldots,a_n])=[b_1,\ldots,b_k,2^{s+1}-1,1,b_k-1,b_{k-1},\ldots,b_1].$$
 By Lemma \ref{nextpar}, $a_n$ satisfies 
 $$
 \begin{cases} [a_1,\ldots,a_{n-1},a_n]<[b_1,\ldots,b_k,2^{s+2}-1,1,b_k-1,b_{k-1},\ldots,b_1], \\ [a_1,\ldots,a_{n-1},a_n+1]>[b_1,\ldots,b_k,2^{s+1}-1,1,b_k-1,b_{k-1},\ldots,b_1] ,\end{cases} 
 $$
and we  need to consider only  two subcases (subcase $n=k$ is impossible since the parity of $n$ and $k$ is different)\\
\begin{equation}\label{final1}
3.1)\,k\le n-1 \,\,\,\,\,\,\,\,\,\,\,\,\,\,\,\,\,\,\,\,\,\,\,\,\,\,\text{and}\,\,\,\,\,\,\,\,\,\,\,\,\,\,\,\,\,\,\, \,\,\,\,\,\,\,\,\,3.2)\,k>n.\end{equation}
 In the case 3.2), analogously to the case  1.2), we come to four options:
\begin{equation}\label{final2}
3.2.1)\,b_n=a_n+1,\,\,\,\,\,\,\,\,\,\,\,
3.2.2)\,b_n=a_n,\,\,\,\,\,\,\,\,\,\,
3.2.3)\,b_n\le a_n-1,\,\,\,\,\,\,\,\,\,\,
3.2.4)\,b_n\ge a_n+2. 
\end{equation}
In every subcase we get a contradiction.\\
In the case 4) we have by Lemma \ref{fold} $$?([a_1,a_2,\ldots,a_n])=[b_1,\ldots,b_k,2^{s+1}-1,1,b_k-1,b_{k-1},\ldots,b_1],$$\\
and by Lemma \ref{nextpar}, $a_n$ satisfies
$$
 \begin{cases} [a_1,\ldots,a_{n-1},a_n]>[b_1,\ldots,b_k,2^{s+2}-1,1,b_k-1,b_{k-1},\ldots,b_1], \\ [a_1,\ldots,a_{n-1},a_n+1]<[b_1,\ldots,b_k,2^{s+1}-1,1,b_k-1,b_{k-1},\ldots,b_1]. \end{cases} 
 $$
Now we consider the same subcases as in \eqref{final1}, and the second subcase splits into the same subcases as in \eqref{final2}.\\
We exhausted all possibilities, getting a contradiction in each of them, hence the assumption \eqref{main} was false. $\square$\\\\
\textbf{Remark 2.} To prove Theorem  \ref{th1} for the greatest fixed point $y$ we should  consider the interval $[y,\frac{1}{2}]$ and the number $\frac{3}{7}$ which belongs to it (by Lemma \ref{2537}) to deduce a contradiction using Lemma \ref{contfunc}.	\\
\textbf{Remark 3.} One can see that a slight generalization of Theorem  \ref{th1} can be proven not only for the smallest and the greatest fixed points, but for any fixed point $x$  which is isolated and unstable at least at one side. \\
For example, $x$ is isolated and unstable at the left side 
if and only if
$$ \exists \varepsilon>0 : \forall y\in (x-\varepsilon,x) \,\,\,\,\text{one has}\,\,\, |y-x|<|?(y)-x| .$$
For the isolated and unstable fixed points instead of Theorem \ref{th1} one can show that the inequality \eqref{eq:2} is valid for all {\it large enough} $n$. 
\section{Proof of Theorem  \ref{th1improved}}
\label{th1improvedproof}
Consider an arbitrary continued fraction $[a_1,\ldots,a_n,\ldots]$. Denote $S_n=a_1+\ldots+a_n$. We need the following lemma from (\cite{Kan}, Theorem 4).
\begin{lem}
\label{kanlem}
Denote $\varphi=\frac{\sqrt{5}+1}{2}$. For any $n\in\mathbb{N}$ one has
\begin{equation}
q_n\le F_{S_n+1}\le\varphi^{S_n}.
\end{equation}
\end{lem}
Now we are ready to prove Theorem \ref{th1improved}. We will only consider the case of the smallest fixed point in the interval $(0,\frac{1}{2})$. The case of the greatest fixed point is treated in the same way.
\begin{solve}
First of all, one can easily see that Theorem \ref{th1improved} holds for $n<36$, as we know the first $36$  partial quotients of $x$. The corresponding sequence is  OEIS A058914 (\cite{OEIS}).
Suppose that $n$ is even. Then by Lemma \ref{contfunc} we have
$$
?\biggl(\frac{p_n}{q_n}\biggr)<\frac{p_n}{q_n}<x=?(x)<?\biggl(\frac{p_{n+1}}{q_{n+1}}\biggr).
$$
Hence
\begin{equation}
\label{mainineqsys}
\frac{1}{(a_{n+1}+1)q_n^2}<x-\frac{p_n}{q_n}<?(x)-?\biggl(\frac{p_n}{q_n}\biggr)<?\biggl(\frac{p_{n+1}}{q_{n+1}}\biggr)-?\biggl(\frac{p_n}{q_n}\biggr)=\frac{1}{2^{S_n+a_{n+1}-1}}.
\end{equation}
We obtain the inequality
\begin{equation}
\label{evencaseineq}
\frac{(a_{n+1}+1)q_n^2}{2^{S_n+a_{n+1}-1}}>1.
\end{equation}
Suppose that $a_{n+1}\ge \kappa_1S_n+2\log_2S_n$. We apply the upper estimate from Lemma \ref{kanlem} and use the fact that $\frac{x}{2^x}$ is strictly decreasing function for $x\ge4$ we obtain
\begin{equation}
1<\frac{(a_{n+1}+1)q_n^2}{2^{S_n+a_{n+1}-1}}<\frac{(\kappa_1 S_n+2\log_2S_n+1)\varphi^{2S_n}}{2^{S_n(\kappa_1+1)+2\log_2S_n-1}}<\frac{2(\kappa_1 S_n+2\log_2S_n+1)}{S_n^2}.
\end{equation}
One can easily see that
$$
\frac{2(\kappa_1 S_n+2\log_2S_n+1)}{S_n^2}<1
$$
for $S_n\ge 4$. We obtain a contradiction.

The case when $n$ is odd is slightly more complicated. Now we have
 $?(\frac{p_{n+1}}{q_{n+1}})<\frac{p_{n+1}}{q_{n+1}}<x=?(x)$. Using the same argument we obtain that
$$
\frac{(a_{n+2}+1)q_{n+1}^2}{2^{S_n+a_{n+1}+a_{n+2}-1}}>1.
$$
As $q_{n+1}<(a_{n+1}+1)q_n$ and $\frac{a_{n+2}+1}{2^{a_{n+2}}}\le1$,
\begin{equation}
\label{oddcaseineqsq}
\frac{(a_{n+1}+1)^2q_{n}^2}{2^{S_n+a_{n+1}-1}}>1.
\end{equation}
Suppose that $a_{n+1}\ge \kappa_1S_n+2\log_2S_n$. Similarly to the previous case, we apply Lemma \ref{kanlem} and use the fact that $\frac{x^2}{2^x}$ is strictly decreasing function for $x\ge 7$. We have
\begin{equation}
\label{prelastineq}
1<\frac{(a_{n+1}+1)^2q_n^2}{2^{S_n+a_{n+1}-1}}<\frac{(\kappa_1 S_n+2\log_2S_n+1)^2\varphi^{2S_n}}{2^{S_n(\kappa_1+1)+2\log_2S_n-1}}<\frac{2(\kappa_1 S_n+2\log_2S_n+1)^2}{S_n^2}.
\end{equation}
From (\ref{prelastineq}) one can easily see that
\begin{equation}
\label{lastineq}
\kappa_1+\frac{2\log_2S_n}{S_n}+\frac{1}{S_n}>\frac{1}{\sqrt{2}}.
\end{equation}
But  (\ref{lastineq}) is not true for $S_n\ge36$ and we obtain a contradiction.
\end{solve}
\noindent
\textbf{Remark 4.} One can prove an even stronger statement, namely
$$
a_{n+1}+\ldots+a_{n+k}< \kappa_1 S_n + 2k\log_2{S_n}.
$$
Combining \eqref{oddcaseineqsq} with $n=m+k-1$ we obtain
\begin{equation}
\label{remark4}
\frac{(a_{m+k}+1)^2q_{m+k-1}^2}{2^{S_{m+k}-1}}>1.
\end{equation}
Let us estimate left-hand side of (\ref{remark4}):
$$
1<\frac{(a_{m+k}+1)^2q_{m+k-1}^2}{2^{S_{m+k}-1}}<\frac{ 2(a_{m+k}+1)^2\cdot\ldots\cdot(a_{m+1}+1)^2q_m^2}{2^{S_m+a_{m+1}+\ldots+a_{m+k}}}.
$$
By assuming that
$
a_{m+1}+\ldots+a_{m+k}\ge \kappa_1 S_m + 2k\log_2{S_m}
$
and applying a similar argument\footnote{We use the fact that the product $\prod\limits_{i=1}^k a_i,\ a_k\in\mathbb{Z}_+$ with the fixed sum of elements $S_k=\sum\limits_{i=1}^k a_i$ attains its maximum when $\bigl|a_i-\frac{S_k}{k}\bigl|<1\ \forall i: 1\le i\le k$.} as in the proof of Theorem \ref{th1improved}, we will get
$$
\frac{\kappa_1}{k}+2\frac{\log_2{S_m}}{S_m}+\frac{1}{S_m}>\frac{1}{2^{1/(2k)}},
$$
and this inequality fails for every $k\ge1$ whenever $S_n\ge36$.\\
\textbf{Remark 5.} Theorem \ref{th1improved} provides a (non-optimal) upper estimate on partial quotients of $x$. However, their mean behavior is much simpler. Denote
\begin{equation}
\begin{split}
\lambda_i=\frac{i+\sqrt{i^2+4}}{2}, \quad\kappa_2=\frac{5\log\lambda_4-4\log\lambda_5}{0.5\log2+\log\lambda_4-\log\lambda_5}\approx4.40104874\ldots.
\end{split}
\end{equation}
In \cite{dushistova} Dushistova, Moshchevitin and Kan proved  
\begin{lem}[\cite{dushistova}, Theorem 3]
\label{MDKtheor}
Let for an irrational number $x$ there exists a constant $C$ such that for all natural $n$ one has
\begin{equation}
\label{kappa2ineq}
S_n\ge\kappa_2n-C.
\end{equation}
Then $?'(x)$ exists and equals $0$. 
\end{lem}
In fact, they showed that if (\ref{kappa2ineq}) holds, then
\begin{equation}
\label{qt2stto2}
\frac{q^2_n}{2^{S_n}}\to 0.
\end{equation}
The fact that $?(x-\delta)<x-\delta$ for any positive $\delta$ implies that $x-\frac{p_n}{q_n}<?(x)-?\bigl(\frac{p_n}{q_n}\bigr)$ for any even $n$. Hence the inequality (\ref{evencaseineq}) holds and we obtain a contradiction with (\ref{qt2stto2}).
 This implies that for any $C$ the inequality $S_n<\kappa_2n+C$ holds for all $n$ large enough. Some calculations show that one can take $C=0$ for all $n\ge 1$. Now we have an obvious consequence of Theorem \ref{th1improved} and Lemma \ref{MDKtheor}.
\begin{foll}
\label{mainfoll}
Let $x$ be fixed point from Theorem \ref{th1}, then 
\begin{equation}
\label{eq:estimate}  
\,\,\, a_{n+1}<\kappa_1\kappa_2n+2\log_2(\kappa_2n)
\end{equation}
for all $ n\ge 2$.
\end{foll}
\section{Proof of Theorem \ref{th2}}
\label{th2proof}
\begin{solve}
From (\ref{evencaseineq}) and (\ref{oddcaseineqsq}) we deduce that for any $n\in\mathbb{N}$ one has
\begin{equation}
\label{allcasesineqsq}
\frac{(a_{n+1}+1)^2q_{n}^2}{2^{S_n+a_{n+1}-1}}>1.
\end{equation}
As $\frac{(a_{n+1}+1)^2}{2^{a_{n+1}-1}}\le\frac{9}{2}$, we have
$$
\frac{2}{9}2^{S_n}<q_n^2
$$
or
$$
S_n<\frac{2}{\log2}\log{q_n}+\log_2{\frac{9}{2}}.
$$
As $a_{n+1}<\kappa_1S_n+2\log_2S_n$, we have 
$$
a_{n+1}<\kappa_1\biggl(\frac{2}{\log2}\log{q_n}+\log_2{\frac{9}{2}}\biggl)+\frac{2}{\log2}\log\biggl(\frac{2}{\log2}\log{q_n}+\log_2{\frac{9}{2}}\biggl).
$$
Consider an arbitrary convergent continued fraction to $x$. As
$$
\biggl|x-\frac{p_n}{q_n}\biggr|>\frac{1}{(a_{n+1}+1)q_n^2},
$$
we see that
$$
\biggl|x-\frac{p_n}{q_n}\biggr|>\frac{1}{\Biggl(\kappa_1\biggl(\frac{2}{\log2}\log{q_n}+\log_2{\frac{9}{2}}\biggl)+\frac{2}{\log2}\log\biggl(\frac{2}{\log2}\log{q_n}+\log_2{\frac{9}{2}}\biggl)+1\Biggr)q_n^2}.
$$
\end{solve}

\section*{Appendix}

\begin{solve}
Suppose that there are at least two fixed points from $(0, \frac{1}{2})$ interval, namely $x_1$ and $x'_1$. Consider an arbitrary rational number $\frac{p}{q}$ between them. Note that the sequence of iterations $z_n=\underbrace{?(?(\ldots?}_{n \text{ iterations}}(\frac{p}{q})\ldots))$ is monotonic and converges to some fixed point in $[x_1, x'_1]$ interval. Suppose that the sequence is decreasing. Denote the limit $\lim\limits_{n\to\infty}z_n$ by $x''_1$, which may coincide with $x_1$. Denote the $n$-th convergent fraction to $x''_1$ by $\frac{p_n}{q_n}$. There exists $N$ such that $\forall n>N$ one has $x''_1<\frac{p_{2n-1}}{q_{2n-1}}<\frac{p}{q}$. As $?(\frac{p}{q})<\frac{p}{q}$, by Lemma \ref{contfunc} we have $x''_1<?(\frac{p_{2n-1}}{q_{2n-1}})<\frac{p_{2n-1}}{q_{2n-1}}$ for all $n>N$. Thus,
\begin{equation}
\label{a2nq2nmin1}
0<\frac{p_{2n-1}}{q_{2n-1}}-?\biggl(\frac{p_{2n-1}}{q_{2n-1}}\biggr)<\frac{p_{2n-1}}{q_{2n-1}}-x''_1<\frac{1}{a_{2n}q_{2n-1}^2}.
\end{equation}

Note that if the sequence $z_n$ is increasing, then using the same argument we obtain that $x''_1>\frac{p_{2n}}{q_{2n}}>\frac{p}{q}$ for all $n$ big enough. The inequality (\ref{a2nq2nmin1}) is replaced by 
\begin{equation}
\label{a2nq2n}
0<?\biggl(\frac{p_{2n}}{q_{2n}}\biggr)-\frac{p_{2n}}{q_{2n}}<x''_1-\frac{p_{2n}}{q_{2n}}<\frac{1}{a_{2n+1}q_{2n}^2}.
\end{equation}

We will finish the proof for decreasing $z_n$ only, because the second case is treated in exactly the same way.  

It follows from (\ref{a2nq2nmin1}) that if $a_{2n}\ge2$ for infinitely many $n>N$, we have a contradiction with (\ref{ge2q2}). Hence there exist $M$ such that $a_{2n}=1$ for $n>M$. We know that
$$
\frac{1}{q^2_{2n-1}}>\frac{p_{2n-1}}{q_{2n-1}}-x''_1>?\biggl(\frac{p_{2n-1}}{q_{2n-1}}\biggr)-?(x''_1).
$$ 
The right-hand side of the previous inequality may be estimated as follows:
$$
?\biggl(\frac{p_{2n-1}}{q_{2n-1}}\biggr)-?(x''_1)>\frac{1}{2^{S_{2n-1}+a_{2n}-1}}-\frac{1}{2^{S_{2n-1}+a_{2n}+a_{2n+1}-1}}>\frac{1}{2^{S_{2n-1}+a_{2n}}}.
$$
Hence, as $a_{2n}=1$ for $n>M$, we have
\begin{equation}
\label{oddcaseineq}
\frac{q_{2n-1}^2}{2^{S_{2n-1}+1}}<1.
\end{equation}
If there exist infinitely many $s$ such that $?\biggl(\frac{p_{2s}}{q_{2s}}\biggl) > \frac{p_{2s}}{q_{2s}}$, then the following holds for $s>N$: \\
$$
\begin{cases} 
0<\frac{p_{2s-1}}{q_{2s-1}}-?\biggl(\frac{p_{2s-1}}{q_{2s-1}}\biggr)<\frac{p_{2s-1}}{q_{2s-1}}-x''_1 \\ 
0<?\biggl(\frac{p_{2s}}{q_{2s}}\biggr)-\frac{p_{2s}}{q_{2s}}<x''_1-\frac{p_{2s}}{q_{2s}}
\end{cases}
$$
A classical theorem states that for any $s\in\mathbb{N}$ at least one of the inequalities $\frac{p_{2s-1}}{q_{2s-1}}-x''_1<\frac{1}{2q^2_{2s-1}}$ or $x''_1-\frac{p_{2s}}{q_{2s}}<\frac{1}{2q^2_{2s}}$ holds. Thus we obtain a contradiction with (\ref{ge2q2}).

Now suppose that for all $n$ big enough we have $?\biggl(\frac{p_{2n}}{q_{2n}}\biggr)<\frac{p_{2n}}{q_{2n}}$. Then, by (\ref{evencaseineq}) and $a_{2n}=1$ we have
\begin{equation}
\label{evencaseineq2}
\frac{(a_{2n+1}+1)q_{2n}^2}{2^{S_{2n-1}+a_{2n+1}}}>1.
\end{equation}
From (\ref{oddcaseineq}) and (\ref{evencaseineq2}) we obtain that
\begin{equation}
\biggl(\frac{q_{2n}}{q_{2n-1}}\biggr)^2\frac{a_{2n+1}+1}{2^{a_{2n+1}-1}}>1.
\end{equation}
Note that
\begin{equation}
\label{almostphi}
\frac{q_{2n}}{q_{2n-1}}=1+[a_{2n-1},a_{2n-2},a_{2n-2},\ldots,a_1]<1+[1,a_{2n-2},1,\ldots,a_1].
\end{equation}
As we already mentioned, $a_{2n}=1$ for all $n>M$. Hence one can estimate the right-hand side of (\ref{almostphi}) as follows
\begin{equation}
\label{almostphi2}
\biggl(\frac{q_{2n}}{q_{2n-1}}\biggr)^2<(\varphi+\varepsilon_n)^2<1.62^2=2.6244 \quad \text{for $n$ large enough}.
\end{equation}
Here $\varepsilon_n$ is some function, which exponentially tends to $0$ as $n$ tends to infinity. Now we can see that $a_{2n+1}\le4$, because if $a_{2n+1}\ge5$, from (\ref{almostphi}) we have
\begin{equation}
\biggl(\frac{q_{2n}}{q_{2n-1}}\biggr)^2\frac{a_{2n+1}+1}{2^{a_{2n+1}-1}}<2.6244\frac{6}{16}<1
\end{equation}
and we obtain a contradiction. Hence there exists $K$ such that for any $n>K$ we have $a_n\le4$. By (\cite{MD}, Theorem 3), we have\footnote{Theorem 3 of that paper states that $?(x)=+\infty$ for numbers $x =[a_1,a_2,\ldots]$ with
$a_i\le 4$ for all $i$, but combining this with Denjoy's formula (\ref{eq:1}) for $?(x)$ it is easily seen that $?(x)=+\infty$ if $\exists N, \forall i > N, a_i\le 4$.} $?'(x''_1)=+\infty$. As 
$$
\frac{?\biggl(\frac{p_{2n-1}}{q_{2n-1}}\biggr)-?(x''_1)}{\frac{p_{2n-1}}{q_{2n-1}}-x''_1}<1
$$
for all $n>N$, we obtain a contradiction. That finishes the proof.
\end{solve}
\noindent
\textbf{Remark 6.} The inequality \eqref{ge2q2} has been verified for all $\frac{p}{q}\in[0,\frac{1}{2}]$ with $q<30000$ and, separately, for 5400 first convergents to the irrational fixed point (or the set of points if Conjecture \ref{mainconj} is false) of $?(x)$ from $(0,\frac{1}{2})$ interval. The only counterexamples, apart from the three trivial rational fixed points $0,\frac{1}{2}$ and $1$, were $\frac{3}{7}$ and $\frac{8}{19}$, both being the convergents to the fixed point. So the numerical calculations allows us to suggest that the inequality \eqref{ge2q2} holds for $q>19$.

\end{document}